\documentclass[reqno]{amsart}
\usepackage{hyperref}
\usepackage{amsmath,amssymb,amsthm,graphicx,mathrsfs,url}
\usepackage[usenames,dvipsnames]{color}
\definecolor{darkred}{rgb}{0.4,0.1,0.1}
\definecolor{darkblue}{rgb}{0.1,0.1,0.4}

\usepackage{tikz}
\usetikzlibrary{datavisualization}
\usetikzlibrary{datavisualization.formats.functions}
\usetikzlibrary{hobby,backgrounds,patterns}

\DeclareMathOperator{\R}{\mathbb R}

\DeclareMathOperator{\cH}{\mathcal H}

\DeclareMathOperator{\dom}{dom}

\def\sa{\mathfrak{a}}

\def\phi{\varphi}

\begin{document}
\title[A variational approach to the hot spots conjecture]
{A variational approach to the hot spots conjecture}

\author[Jonathan Rohleder \hfil \hfilneg] {Jonathan Rohleder}  

\address{Jonathan Rohleder \newline
Matematiska institutionen \\ Stockholms universitet \\
106 91 Stockholm \\
Sweden}
\email{jonathan.rohleder@math.su.se}


\subjclass[2000]{} \keywords{Laplacian, eigenfunctions, hot spots conjecture, spectral theory}
\begin{abstract}
We review a recent new approach to the study of critical points of Laplacian eigenfunctions. Its core novelty is a non-standard variational principle for the eigenvalues of the Laplacians with Neumann and Dirichlet boundary conditions on bounded, simply connected planar domains. This principle can be used to provide simple proofs of some previously known results on the hot spots conjecture.
\end{abstract}

\maketitle \numberwithin{equation}{section}
\newtheorem{theorem}{Theorem}[section]
\newtheorem{corollary}[theorem]{Corollary}
\newtheorem{proposition}[theorem]{Proposition}
\newtheorem{lemma}[theorem]{Lemma}
\newtheorem{remark}[theorem]{Remark}
\newtheorem{problem}[theorem]{Problem}
\newtheorem{example}[theorem]{Example}
\newtheorem{definition}[theorem]{Definition}
\allowdisplaybreaks

\section{Introduction}

The {\it hot spots conjecture} seems to have been formulated first by Rauch in 1974, see Section~II.5 in Kawohl's book \cite{K85}. It suggests, in its strongest form, that any eigenfunction corresponding to the smallest positive eigenvalue of the Laplacian on a bounded domain with Neumann boundary conditions attains its maximum and minimum only on the boundary. This would imply that the hottest and coldest spots in an insulated body with a ``generic'' initial heat distribution diverge from each other for large time as much as they can, i.e.\ converge to the boundary; cf.~\cite{K85}.

Let us briefly review the major steps forward on this conjecture. First of all, the conjecture is true for simple shapes on which the eigenfunctions in question can be computed explicitly, such as balls, cubes, equilateral triangles or right isosceles triangles. Moreover, there exist by now a few non-trivial results for domains in the plane. Ba\~nuelos and Burdzy \cite{BB99} proved that the conjecture is true for obtuse triangles and sufficiently long convex domains with symmetries, and Atar and Burdzy~\cite{AB04} showed it for so-called lip-domains, i.e.\ domains enclosed by the graphs of two Lipschitz continuous functions with Lipschitz constants at most one. Jerison and Nadirashvili proved it for domains symmetric w.r.t.\ both coordinate axes \cite{JN00} for which all horizontal and vertical cross sections are intervals. On the other hand Burdzy and Werner \cite{BW99} and Burdzy \cite{B05} constructed counterexamples given by certain multiply-connected domains, see also the numerical study \cite{K21}. The most recent advances include a prove for general triangles by Judge and Mondal \cite{JM20,JM20err}, which was preceeded by Siudeja's partial result \cite{S15}. For further recent related results we refer the reader to, e.g., \cite{JM22,KT19,MPW23,S20,S23}. The conjecture is still open for general simply connected domains in the plane, as well as in higher dimensions, and is believed to be true at least for convex domains.

It seems that most proofs of the key results on the hot spots conjecture are either based on probabilistic methods, exploiting reflected Brownian motion \cite{AB04,BB99,BB00,B05,BW99,P02,S20}, or rely on tracing critical points of eigenfunctions under perturbations of the domain \cite{JN00,JM20}. The author of this note is suggesting a completely independent approach \cite{R23+}, inspired by the following classical observation: as is well known, in order to study nodal sets of eigenfunctions, it is convenient to make use of the fact that these eigenfunctions are optimizers of certain variational principles. For instance, the first non-trivial eigenvalue $\mu_2$ of the Neumann Laplacian $- \Delta_{\rm N}$ on a bounded Lipschitz domain $\Omega$ is given by
\begin{align}\label{eq:minMaxclassic}
 \mu_2 = \min_{\substack{\psi \in H^1 (\Omega) \setminus \{0\} \\ \int_\Omega \psi = 0}} \frac{\int_\Omega |\nabla \psi|^2}{\int_\Omega |\psi|^2},
\end{align}
and a function $\psi$ in the Sobolev space $H^1 (\Omega)$ with vanishing integral is a minimizer of \eqref{eq:minMaxclassic} if and only if $\psi$ is an eigenfunction of $- \Delta_{\rm N}$ corresponding to $\mu_2$. From this characterization it can be derived easily that $\psi$ has precisely two nodal domains; see, e.g.,~\cite[Chapter VI, § 6]{CH}. Therefore in order to study critical points of $\psi$, it seems natural to search for variational principles in which the gradient of $\psi$, instead of $\psi$ itself,  is a minimizer. The following theorem was obtained in \cite{R23+}.

\begin{theorem}\label{thm:minMaxIntro}
Assume that $\Omega \subset \R^2$ is a bounded, simply connected Lipschitz domain with piecewise $C^\infty$-smooth boundary and that all its corners, if any, are convex. Then
\begin{align}\label{eq:minPrincipleNeumannIntro}
 \mu_2 = \min_{u = \binom{u_1}{u_2} \in \cH_{\rm N} \setminus \{0\}} \frac{\int_\Omega \left( |\nabla u_1|^2 + |\nabla u_2|^2 \right) - \int_{\partial \Omega} \kappa  \big( |u_1|^2 + |u_2|^2 \big)}{\int_\Omega \left( |u_1|^2 + |u_2|^2 \right)},
\end{align}
where $\cH_{\rm N}$ consists of all vector fields with components in the Sobolev space $H^1 (\Omega)$ such that their traces satisfy $\langle u, \nu \rangle = 0$ a.e.\ on $\partial \Omega$, and $\kappa$ is the signed curvature on $\partial \Omega$ w.r.t.\ the outer unit normal $\nu$, defined on all boundary points except corners. The minimizers of \eqref{eq:minPrincipleNeumannIntro} are precisely the gradients of eigenfunctions $\psi$ of the Neumann Laplacian corresponding to $\mu_2$.
\end{theorem}

In the following Section \ref{sec:var} we will sketch how this variational principle is obtained; it turns out to be a particular case of a min-max principle for all eigenvalues of the Laplacians with Neumann and Dirichlet boundary conditions. In the final Section~\ref{sec:hotSpots}, we explain how several previously known results on the hot spots conjecture can be derived in an elementary fashion from Theorem \ref{thm:minMaxIntro}.

\section{A non-standard variational principle}\label{sec:var}

In this section we sketch the construction which leads to Theorem \ref{thm:minMaxIntro}. The main idea is to construct a self-adjoint operator, acting as the negative Laplacian on vector fields, for which $\nabla \psi$ is an eigenfunction if $\psi$ is a non-constant eigenfunction of the Neumann Laplacian. 

Let us briefly fix some notation; for more details we refer the reader to \cite{R23+}. On a bounded, connected Lipschitz domain $\Omega \subset \R^2$ we denote by $- \Delta_{\rm N}$ the Neumann Laplacian, i.e.\
\begin{align*}
 - \Delta_{\rm N} u & = - \Delta u, \quad \dom (- \Delta_{\rm N}) = \left\{u \in H^1 (\Omega) : \Delta u \in L^2 (\Omega), \partial_\nu u |_{\partial \Omega} = 0 \right\},
\end{align*}
and by $- \Delta_{\rm D}$ the Dirichlet Laplacian,
\begin{align*}
 - \Delta_{\rm D} u & = - \Delta u, \quad \dom (- \Delta_{\rm D}) = \left\{u \in H^1 (\Omega) : \Delta u \in L^2 (\Omega), u |_{\partial \Omega} = 0 \right\}.
\end{align*}
The boundary conditions of the operators have to be understood in an appropriate weak sense; $u |_{\partial \Omega}$ denotes the trace of $u$ on $\partial \Omega$ and $\partial_\nu u |_{\partial \Omega}$ is the derivative of $u$ on the boundary in the direction of the outward pointing normal vector $\nu$. Both operators are unbounded and self-adjoint in $L^2 (\Omega)$, and their spectra consist of isolated eigenvalues of finite multiplicities. Let
\begin{align}\label{eq:EVN}
 0 = \mu_1 < \mu_2 \leq \mu_3 \leq \dots 
\end{align}
be an enumeration of the eigenvalues of $- \Delta_{\rm N}$ and let
\begin{align}\label{eq:EVD}
 \lambda_1 < \lambda_2 \leq \lambda_3 \leq \dots 
\end{align}
be the eigenvalues of $- \Delta_{\rm D}$, both counted according to their multiplicities. The eigenfunctions of $- \Delta_{\rm N}$ corresponding to $\mu_1 = 0$ are the constant functions. 

Let us now assume, in addition, that $\Omega$ has a piecewise smooth boundary. In the space $L^2 (\Omega)^2$ of square-integrable two-component vector fields on $\Omega$ we define the sesquilinear form
\begin{align*}
 \sa \left[ u, v \right] & = \int_\Omega \big( \left\langle\nabla u_1, \nabla v_1\right\rangle + \left\langle\nabla u_2, \nabla v_2\right\rangle \big) - \int_{\partial \Omega} \kappa \langle u, v \rangle, \quad u = \binom{u_1}{u_2}, v =  \binom{v_1}{v_2},
\end{align*}
with domain 
\begin{align*}
 \dom \sa := \cH_{\rm N} := \left\{ u \in H^1 (\Omega)^2 : \left\langle u |_{\partial \Omega}, \nu \right\rangle = 0~\text{on}~\partial \Omega \right\}.
\end{align*}
Here we denote by $\kappa$ the signed curvature function along the piecewise smooth curve $\partial \Omega$ w.r.t.\ the outer unit normal $\nu$, defined on all points of $\partial \Omega$ except possible corners; in particular, $\kappa$ is bounded and piecewise smooth, with possible jumps at the corners. If $\Omega$ is convex then $\kappa (x) \leq 0$ holds for almost all $x \in \partial \Omega$. 

The sesquilinear form $\sa$ is symmetric, semi-bounded, densely defined in $L^2 (\Omega)^2$ and closed. Hence, by \cite[VI, Theorem~2.1]{Kato} there exists a unique self-adjoint operator $A$ in $L^2 (\Omega)^2$ such that $\dom A \subset \dom \sa$ and
\begin{align*}
 (A u, v)_{L^2 (\Omega)^2} = \sa [u, v], \quad u \in \dom A, v \in \dom \sa,
\end{align*}
where $(\cdot, \cdot)_{L^2 (\Omega)^2}$ stands for the standard inner product in the space $L^2 (\Omega)^2$. It can be computed that that the operator $A$ acts as the Laplacian,
\begin{align*}
 A u = \binom{- \Delta u_1}{- \Delta u_2}, \quad u \in \dom A,
\end{align*}
and its domain consists of sufficiently regular vector fields satisfying the conditions
\begin{align*}
 \langle u |_{\partial \Omega}, \nu \rangle = 0 \quad \text{and} \quad \partial_1 u_2 - \partial_2 u_1 = 0 \quad \text{on}~\partial \Omega,
\end{align*}
interpreted in a weak sense; cf.\ \cite[Lemma~3.3 and Remark~3.4]{R23+}. The operator $A$ is intimately related to the operators $- \Delta_{\rm N}$ and $- \Delta_{\rm D}$ in the following way. We are imposing here conditions on the domain $\Omega$ which make sure that functions in the domains of $- \Delta_{\rm D}$ and $- \Delta_{\rm N}$ belong to the Sobolev space $H^2 (\Omega)$. We make use of the notation $\nabla^\perp \phi = (- \partial_2 \phi, \partial_1 \phi)^\top$.

\begin{theorem}\label{thm:translateEV}
Assume that $\Omega \subset \R^2$ is a bounded Lipschitz domain with piecewise $C^\infty$-smooth boundary whose corners are convex. Let $A$ be the self-adjoint operator in $L^2 (\Omega)^2$ associated with the sesquilinear form $\sa$. Moreover, let the eigenvalues of $- \Delta_{\rm N}$ be enumerated as in \eqref{eq:EVN} and let $\psi_1, \psi_2, \dots$ form an orthonormal basis of $L^2 (\Omega)$ such that $- \Delta_{\rm N} \psi_k = \mu_k \psi_k$ holds for $k = 1, 2, \dots$; analogously let the eigenvalues of $- \Delta_{\rm D}$ be enumerated as in \eqref{eq:EVD} and let $\phi_1, \phi_2, \dots$ be an orthonormal basis of $L^2 (\Omega)$ consisting of corresponding eigenfunctions, $- \Delta_{\rm D} \phi_k = \lambda_k \phi_k$ for all $k$. Then the following hold.
\begin{enumerate}
 \item For each $k \geq 2$, $\nabla \psi_k$ is nontrivial, belongs to $\dom A$, and satisfies $A \nabla \psi_k = \mu_k \nabla \psi_k$. Moreover, the functions $\frac{1}{\sqrt{\mu_k}} \nabla \psi_k$ form an orthonormal basis of $\nabla H^1 (\Omega)$.
 \item For each $k \geq 1$, $\nabla^\perp \phi_k$ is nontrivial, belongs to $\dom A$, and satisfies $A \nabla^\perp \phi_k = \lambda_k \nabla^\perp \phi_k$. Moreover, the functions $\frac{1}{\sqrt{\lambda_k}} \nabla^\perp \phi_k$ form an orthonormal basis of $\nabla^\perp H_0^1 (\Omega)$.
\end{enumerate}
In particular, if $\Omega$ is simply connected then the spectrum of $A$ coincides with the union of the positive eigenvalues of the Neumann and Dirichlet Laplacians, counted with multiplicities.
\end{theorem}

For a rigorous proof we refer the reader to \cite[Theorem~3.2]{R23+}. For a sketch, let $\psi$ be any non-constant eigenfunction of $- \Delta_{\rm N}$ corresponding to an eigenvalue $\mu$ and $u = \nabla \psi$. Then $u$ belongs to $\dom A$: for the boundary conditions, note that
\begin{align*}
 u \cdot \nu = \nabla \psi \cdot \nu = 0 \quad \text{on}~\partial \Omega,
\end{align*}
and
\begin{align*}
 \partial_1 u_2 - \partial_2 u_1 = \partial_1 \partial_2 \psi - \partial_2 \partial_1 \psi = 0
\end{align*}
by the Schwartz lemma, and this holds in particular on $\partial \Omega$. Moreover, it is clear that
\begin{align*}
 A u = \binom{- \Delta \partial_1 \psi}{- \Delta \partial_2 \psi} = - \nabla \Delta \psi = \mu u.
\end{align*}
Similarly, if $- \Delta_{\rm D} \phi = \lambda \phi$ and $u = \nabla^\perp \phi$, then $u$ satisfies the boundary conditions required in $\dom A$: firstly,
\begin{align*}
 u \cdot \nu = \nabla^\perp \phi \cdot \nu = 0 \quad \text{on}~\partial \Omega
\end{align*}
as this corresponds to the tangential derivative on $\partial \Omega$ and $\phi = 0$ constantly there. Secondly,
\begin{align*}
 \partial_1 u_2 - \partial_2 u_1 = \partial_1 \partial_1 \phi + \partial_2 \partial_2 \phi = \Delta \phi = - \lambda \phi = 0 \quad \text{on}~\partial \Omega,
\end{align*}
making use of the differential equation and the boundary condition for $\phi$. It follows $A \nabla^\perp \phi = \lambda \nabla^\perp \phi$. Furthermore, proof of the orthonormal basis properties mentioned in the theorem is straightforward.

Let us now assume, in addition, that $\Omega$ is simply connected. Then the famous Helmholtz decomposition reads
\begin{align*}
 L^2 (\Omega)^2 = \nabla H^1 (\Omega) \oplus \nabla^\perp H_0^1 (\Omega),
\end{align*}
see, e.g., \cite[Lemma 2.10]{K10}. Using this, it follows from Theorem \ref{thm:translateEV} (1) and (2) that the spectrum of $A$ consists precisely of all non-trivial eigenvalues of $- \Delta_{\rm N}$ and $- \Delta_{\rm D}$, including multiplicities. As it is well-known that $\mu_2 < \lambda_1$ (see \cite{F05,P52}), we obtain the following min-max principle.

\begin{theorem}\label{thm:minMaxDN}
Assume that $\Omega \subset \R^2$ is a bounded, simply connected domain with piecewise $C^\infty$-smooth boundary and that all corners, if any, are convex. Denote by
\begin{align*}
 \eta_1 \leq \eta_2 \leq \dots
\end{align*}
the union of the positive eigenvalues of $- \Delta_{\rm N}$ and $- \Delta_{\rm D}$, counted according to their multiplicities. Then
\begin{align*}
 \eta_j = \min_{\substack{F \subset \cH_{\rm N}~\text{subspace} \\ \dim F = j}} \max_{u = \binom{u_1}{u_2} \in F \setminus \{0\}} \frac{\int_\Omega \left( |\nabla u_1|^2 + |\nabla u_2|^2 \right) - \int_{\partial \Omega} \kappa  \big( |u_1|^2 + |u_2|^2 \big)}{\int_\Omega \left( |u_1|^2 + |u_2|^2 \right)}.
\end{align*}
Especially, the first positive eigenvalue $\mu_2$ of $- \Delta_{\rm N}$ is given by
\begin{align}\label{eq:minPrincipleNeumann}
 \mu_2 = \min_{u = \binom{u_1}{u_2} \in \cH_{\rm N} \setminus \{0\}} \frac{\int_\Omega \left( |\nabla u_1|^2 + |\nabla u_2|^2 \right) - \int_{\partial \Omega} \kappa  \big( |u_1|^2 + |u_2|^2 \big)}{\int_\Omega \left( |u_1|^2 + |u_2|^2 \right)}.
\end{align}
Moreover, if $\psi$ is an eigenfunction of $- \Delta_{\rm N}$ corresponding to $\mu_2$ then the minimum in \eqref{eq:minPrincipleNeumann} is attained at $u = \nabla \psi$, and, conversely, each minimizer $u$ of \eqref{eq:minPrincipleNeumann} satisfies $u = \nabla \psi$ for some $\psi \in \ker (- \Delta_{\rm N} - \mu_2)$.
\end{theorem}

\section{Application to the hot spots conjecture}\label{sec:hotSpots}

In this section we indicate how the variational principle of Theorem \ref{thm:minMaxDN} can be applied to the hot spots conjecture. The first result discussed here was first proven by Atar and Burdzy \cite{AB04} by probabilistic methods. It comprises the class of so-called lip domains introduced by Burdzy and Chen in \cite{BC98}.

\begin{definition}\label{def:lip}
A bounded Lipschitz domain $\Omega \subset \R^2$ is called {\em lip domain} if
\begin{align*}
 \Omega = \left\{ (x, y)^\top: f_1 (x) < y < f_2 (x), x \in (a, b) \right\},
\end{align*}
where $f_1, f_2 : [a, b] \to \R$ are Lipschitz continuous functions with Lipschitz constant at most one such that $f_1 (x) < f_2 (x)$ for all $x \in (a, b)$, $f_1 (a) = f_2 (a)$ and $f_1 (b) = f_2 (b)$.
\end{definition}

\begin{figure}[h]
\begin{center}
\begin{tikzpicture}[rotate=-45,transform shape]
\coordinate (a) at (6,-1);
\coordinate (b) at (10, -0.5);
\coordinate (c) at (11,-0.1);
\coordinate (d) at (11,3.24);
\coordinate (e) at (11,1);
\draw[thick] (a) edge (b);
\draw[thick] (e) edge (d);
\node[rotate=45] at (8.5,0.6) {$\Omega$};
\draw[scale=1,domain=6:11,smooth,variable=\x,thick] plot ({\x},{sqrt(\x - 6) + 1});
\draw[scale=1,domain=10:11,smooth,variable=\x,thick] plot ({\x},{1.5*(sqrt(\x - 10) - 0.5/1.5)});
\draw[thick] (a) edge (6,1);
\end{tikzpicture}
\end{center}
\caption{A lip domain.}
\label{fig:domains}
\end{figure}
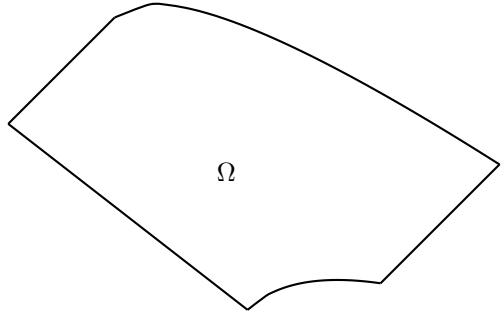

Figure~\ref{fig:domains} shows a typical lip domain. Other examples include right and obtuse triangles (in contrast to acute triangles) or right and obtuse trapezoids (in contrast to acute trapezoids). Clearly, lip domains are simply connected.

\begin{theorem}\label{thm:main}
Assume that $\Omega$ is a lip domain with piecewise $C^\infty$-smooth boundary whose corners are convex. Then the following assertions hold.
\begin{enumerate}
 \item If $\Omega$ is not a square then the first positive eigenvalue $\mu_2$ of the Neumann Laplacian on $\Omega$ is simple, i.e.\ the corresponding eigenfunction $\psi$ is unique up to scalar multiples.
 \item If $\Omega$ is not a rectangle then $\psi$ may be chosen such that its directional derivatives in both directions $\mathbf{e}_1 + \mathbf{e}_2$ and $\mathbf{e}_1 - \mathbf{e}_2$ are positive inside $\Omega$, where $\mathbf{e}_1$ and $\mathbf{e}_2$ are the standard basis vectors in $\R^2$. In particular, $\psi$ does not have any critical point inside $\Omega$ and, hence, takes its maximum and minimum on $\partial \Omega$ only.
\end{enumerate}
\end{theorem}

The proof relies on the fact that rotating any lip domain by $\pi/4$ in positive direction leads to a domain for which the outer unit normal vector field $\nu$ on the boundary satisfies 
\begin{align}\label{eq:oppositeQuadrants}
 \nu (x) \in \overline{Q_2} \cup \overline{Q_4},
\end{align}
that is, $\nu (x)$ belongs to the union of the closed second and fourth quadrants in the plane, for almost all $x \in \partial \Omega$. We will now sketch the proof of Theorem \ref{thm:main} assuming that $\Omega$ is rotated such that \eqref{eq:oppositeQuadrants} holds for almost all $x \in \partial \Omega$.

Let $\psi$ be any eigenfunction of $- \Delta_{\rm N}$ corresponding to $\mu_2$. Then $u = \nabla \psi$ minimizes \eqref{eq:minPrincipleNeumann}. Take $v = (|u_1|, |u_2|)^\top$. Then 
\begin{align*}
 \|v\|_{L^2 (\Omega)^2} = \|u\|_{L^2 (\Omega)^2} \quad \text{and} \quad \sa [v] \leq \sa [u].
\end{align*}
Moreover, the condition \eqref{eq:oppositeQuadrants} guarantees that $v |_{\partial \Omega} \cdot \nu = 0$ holds, since the components of $\nu (x)$ have opposite signs for almost all $x \in \partial \Omega$. Hence, $v \in \cH_{\rm N}$ and $v$ is another minimizer of \eqref{eq:minPrincipleNeumann}, whose components are non-negative everywhere in $\Omega$. Thus $v \in \ker (A - \mu_2)$ and there exists an eigenfunction $\psi'$ of $- \Delta_{\rm N}$ corresponding to the eigenvalue $\mu_2$ such that $\nabla \psi' = v$. However,
\begin{align*}
 \Delta \partial_j \psi' = - \mu_2 |u_j| \leq 0, \quad j = 1, 2,
\end{align*}
and, thus, the minimum principle for superharmonic functions yields that $v_j = \partial_j \psi'$ is either constantly zero or strictly positive in $\Omega$, $j = 1, 2$. As $\psi'$ being constant in one direction is only possible if $\Omega$ is a rectangle, the assertion (2) follows. For assertion (1), if $\psi$ and $\psi'$ are linearly independent and $\Omega$ is not a square, one can construct a linear combination of $\nabla \psi$ and $\nabla \psi'$ vanishing at some point in $\Omega$, contradicting the reasoning in the proof of assertion (2). A complete proof can be found in \cite[Section~4]{R23+}.

Theorem \ref{thm:minMaxDN} can be employed for a more careful analysis of the nodal lines of the components of $\nabla \psi$, for an eigenfunction $\psi$ of $- \Delta_{\rm N}$. This can be used, for instance, to reprove the following theorem of Jerison and Nadirashvili \cite[Theorem~1.1]{JN00}.

\begin{theorem}\label{thm:JN}
Assume that $\Omega \subset \R^2$ has a piecewise $C^\infty$-smooth boundary without non-convex corners and is symmetric with respect to both coordinate axes. Furthermore, assume that all vertical and horizontal cross sections of $\Omega$ are intervals and that $\Omega$ is not a rectangle. Then the following hold.
\begin{enumerate}
 \item For any eigenfunction $\psi$ corresponding to $\mu_2$ that is odd w.r.t.\ $x$ and even w.r.t.\ $y$, $\partial_x \psi$ does not have any zero in $\Omega$; if $\psi$ is chosen such that $\psi > 0$ if $x > 0$, then $\partial_x \psi > 0$ in $\Omega$. Moreover, $\partial_y \psi$ vanishes exactly on the axes.
 \item For any eigenfunction $\psi$ corresponding to $\mu_2$ that is even w.r.t.\ $x$ and odd w.r.t.\ $y$, $\partial_y \psi$ does not have any zero in $\Omega$; if $\psi$ is chosen such that $\psi > 0$ if $y > 0$, then $\partial_y \psi > 0$ in $\Omega$. Moreover, $\partial_x \psi$ vanishes exactly on the axes.
\end{enumerate}
\end{theorem}

A proof of Theorem \ref{thm:JN} based on Theorem \ref{thm:minMaxDN} will be published elsewhere.

It should, finally, be pointed out that the regularity assumptions on $\partial \Omega$ in the above theorems are stronger than those in \cite{AB04,JN00}.

\section*{Acknowledgements}
This research is financially supported by grant no.\ 2022-03342 of the Swedish Research Council (VR).

\end{document}